\documentclass[11pt]{article}
\usepackage{amssymb,amsmath,amsthm,latexsym, bbm}
\usepackage{enumerate}
\usepackage{color}
\usepackage{graphicx}
\usepackage{mathrsfs}
\usepackage[mathscr]{euscript}
\usepackage{amscd}

\newtheorem{thm}{Theorem}[section]

\newtheorem{lemma}[thm]{Lemma}

\newtheorem{defn}[thm]{Definition}

\numberwithin{equation}{section}

\def\Ex {{\mathbb E}}

\def\R {{\mathbb R}}

\def\P {{\mathbb P}}

\def\K{{\mathbf K}}
\def\J{{\mathbf J}}

\newcommand{\norm}[1]{\lVert#1\rVert}

\def\wt{\widetilde}
\def\wh{\widehat}

\def\op{\mathcal L}
\def\1{{\mathbf 1}}

\def\<{{\langle}}
\def\>{{\rangle}}
\def\eps{\varepsilon}
\def\wh{\widehat}

\def\proof{{\medskip\noindent {\bf Proof. }}}

\def\qed{{\hfill $\square$ \bigskip}}

\def\Levy {{L\'evy\ }}

\begin{document}
\title{Stability of Heat Kernel Estimates for Diffusions with Jumps under Non-local Feynman-Kac Perturbations}

\author{{Zhen-Qing Chen} \quad \hbox{and} \quad {Lidan Wang}}

\date{\today }

\maketitle

\begin{abstract}
In this paper we show that the two-sided 
 heat kernel estimates for a class
of (not necessarily symmetric) diffusions with jumps   are stable under non-local
Feynman-Kac perturbations.  
\end{abstract}

	

\bigskip

\section{Introduction}

Suppose that $X$ is a Hunt process on a state space $E$ with transition semigroup
$\{P_t: t\ge 0\}$. A Feynman-Kac transform of $X$ is given by
\begin{equation}\label{e:1.1}
 T_t f (x) =\Ex_x \left[ {\rm Exp} (C_t) f(X_t)\right],
\end{equation}
where $C_t$ is an additive functional of $X$ of finite variation. 
When $C_t$ is a continuous additive functional
of $X$, the transform above is called a (local) Feynman-Kac transform. 
Feynman-Kac transforms play an important role in the probabilistic
as well as analytic aspect of potential theory, and also in mathematical physics. 
For example, Feynman-Kac transforms for Brownian motion on Euclidean spaces have been
studied extensively, see \cite{BM,  CZ, S} and the references therein for a survey on this topic.
When $X$ is discontinuous, additive functionals of $X$ of finite variation can be discontinuous
and there are many of them. When $C_t$ is a discontinuous additive functional of $X$, the transform of \eqref{e:1.1} is called a non-local Feynman-Kac transform.
Non-local Feynman-Kac transforms have received quite lot of attention recently in connection with the study of potential theory for discontinuous Markov processes and non-local operators; see, for example,
\cite{C1, C2, C3, CFKZ, CKi, CKS, CS1, CS2, W} and the references therein. 
An important question related to Feynman-Kac transforms is the stability
of various properties. In particular, stability of heat kernel estimates for purely discontinuous Markov processes under non-local Feynman-Kac perturbations have been studied in \cite{CKS, W}.  See  \cite{BS} for a related work.  
In this paper, we study the stability of heat kernel estimates for diffusions with jumps under non-local Feynman-Kac perturbations.  

A generic strong Markov process may have both the continuous (diffusive) part and the purely discontinuous (jumping) part. In Chen and Kumagai \cite{CKS}, symmetric diffusion processes with jumps on $\R^d$
having generators
\begin{equation}
	\label{generator1}
	\op  u(x)=\frac{1}{2}\sum_{i,j=1}^d\frac{\partial}{\partial x_i}\Big(a_{ij}(x)\frac{\partial u(x)}{\partial x_j}\Big)+\lim_{\eps\downarrow0}\int_{\{y\in\R^d:|y-x|>\eps\}}\big(u(y)-u(x)\big)\frac{c(x,y)}{|x-y|^{d+\alpha}}dy
\end{equation}
are studied, 
where $\alpha\in(0,2)$, $A(x)=(a_{ij}(x))_{1\leq i,j\leq d}$ is a measurable $d\times d$ matrix-valued function on $\R^d$ that is uniformly elliptic and bounded, and $c(x,y)$ is symmetric function that is bounded between two positive constants. It is shown that there is a Feller process $X$ having strong Feller property associated with $\op$, which we call symmetrid diffusion with jumps. The Feller process 
$X_t$ has a jointly H\"older continuous transition density function $p(t, x, y)$ with respect to the Lebesgue measure on $\R^d$ and the following two-sided estimates hold. 
There exist positive constants $c_k$, $1\leq k\leq 4$,  such that for every $t>0$ and $x,y\in\R^d$,
\begin{eqnarray} \label{e:1.3}
&& c_1\left(t^{-d/2}\wedge t^{-d/\alpha}\right)\wedge\left(t^{-d/2}{\rm Exp}\left(-\frac{c_2|x-y|^2}{t}\right)+t^{-d/\alpha}\wedge \frac{t}{|x-y|^{d+\alpha}}\right)\nonumber \\  
  &\leq  & p(t,x,y) 
\leq c_3\left(t^{-d/2}\wedge t^{-d/\alpha}\right)\wedge\left(t^{-d/2}{\rm Exp}\left(-\frac{c_4|x-y|^2}{t}\right)+t^{-d/\alpha}\wedge \frac{t}{|x-y|^{d+\alpha}}\right)  .  
\end{eqnarray} 
For $a, b\in \R$, $a\wedge b:=\min\{a, b\}$ and  $a\vee b:=\max \{a, b\}$.  

\medskip

Recently the following non-symmetric non-local operator  
\begin{equation}
	\label{generator2}
	\op u(x)=\frac{1}{2}\sum_{i,j=1}^{d}a_{ij}(x)\frac{\partial^2 u(x)}{\partial x_i \partial x_j}
	+\sum_{j=1}^d b_j (x) \frac{\partial u(x)}{\partial x_j} + 
	\int_{\R^d}\left(u(x+z)-u(x)-1_{\{|z|\leq 1\}}\cdot\nabla u(x)\right)\frac{c(x,z)}{|z|^{d+\alpha}}dz,
\end{equation}
has been studied in \cite{CHXZ}, 
where $\alpha\in(0,2)$, and $A=(a_{ij}(x))_{1\leq i,j\leq d}$ is a measurable $d\times d$ matrix-valued function on $\R^d$
that is uniformly elliptic  and bounded, 
and is H\"older continuous , $b(x)=(b_1(x), \dots, b_d (x))$ is an $\R^d$-valued function that is in some Kato class,
and $c(x,z)\geq 0$ is a  bounded measurable function such that when $\alpha =1$, 
$$
\int_{\{r<|z|<R\}} z \, c(x, z) dz =0 \quad \hbox{for every } 0<r<R<\infty .
$$
Chen, Hu, Xie and Zhang \cite{CHXZ} showed, among other things, that there is a Feller process $X_t$ 
having strong Feller property associated with the above generator, and $X_t$ has a jointly continuous transition density function $p(t, x, y)$ with respect to the Lebesgue measure on $\R^d$. Moreover, when $c(x, z)$ is bounded between two positive constants, the two-sided estimates \eqref{e:1.3} are shown to 
hold for $p(t, x, y)$ on $(0, T]\times \R^d \times \R^d$
for every $T>0$.  In fact, more general time-dependent operators of the form \eqref{generator2} 
are studied in \cite{CHXZ}. 

In this paper, we start with a Hunt process $X$ on $\R^d$  with $d\geq 2$
that has a jointly continuous transition density function $p(t, x, y)$ that
enjoys two-sided estimates \eqref{e:1.3} on $(0, T]\times \R^d \times \R^d$. 
Under this assumption,  the Hunt process $X$ has a L\'evy system $(N(x, dy), dt)$ with
$N(x, dy)=\frac{c(x, y)}{|x-y|^{d+\alpha}} dy$ for some measurable function $c(x, y)$ bounded between 
two positive constants; see \eqref{e:3.7} below. That is,
for every non-negative function $\varphi (x, y)$ on $\R^d \times \R^d$ that vanishes along the diagonal,
$$
\Ex_x \Big[ \sum_{0<s\leq t} \varphi (X_{s-}, X_s) \Big] = \Ex_x \left[  \int_0^t \int_{\R^d} \varphi (X_s, y) N(X_s, dy) ds
\right], \quad x\in \R^d, \ t>0.
$$
Here we use the convention that we extend the definition of functions to cemetery point $\partial$ by setting 0 value there;
for example $\varphi (x, \partial )=0$. 
For convenience, we take $T=1$. 
We will study the stability of heat kernel estimates under non-local Feynman-Kac transform:  
$$
T_tf(x)=\Ex_x\Big[\exp\Big(A^\mu_t+\sum_{s\leq t}F(X_{s-},X_s)\Big)f(X_t)\Big],
$$
where $A^\mu$ is a continuous additive functional of $X$ of finite variations having signed Revuz measure $\mu$
and  $F(x,y)$ is a bounded measurable function vanishing on the diagonals. We point out that in this paper 
we do not require $F$ to be symmetric. 
Informally, the semigroup $(T_t^{\mu,F};t\geq0)$ has  generator
\begin{equation}
	\label{generatorFK}
	\mathcal Af(x)=(\op+\mu)f(x)+  \int_{\R^d}\left(e^{F(x,y)}-1\right)  f(y) N(x,  dy) ,
\end{equation}
where $\op$ is the   infinitesimal generator of $X$;  
see \cite[Remark 1 on p.53]{CS1} for a calculation. 
We show that if $\mu$ and $F$ are in certain Kato class of $X$,
the non-local Feynman-Kac semigroup $\{T_t; t\geq 0\}$ has a heat kernel $q(t, x, y)$
and $q(t, x, y)$ has two-sided estimates \eqref{e:1.3} on $(0, 1]\times \R^d \times \R^d$
but with a set of possible different constants $c_k$, $1\leq k\leq 4$.
Comparing with \cite{CKS, W, BS}, the novelty of this paper is that $X$ has both the diffusive and jumping  components,
and that the Gaussian bounds in \eqref{e:1.3} have different constants $c_2$ and $c_4$ in the exponents
for the upper and lower bound estimates.  These features made the perturbation estimates more challenging.

The rest of the paper is organized as follows. 
Section \ref{sec:preliminary} gives the basic setup of the problem and the statement of the main results of this paper.
In Section \ref{sec:3P}, we various 3P type inequalities needed to study non-local Feynman-Kac perturbations.
Proof of the main results,  the two-sided estimates for the heat kernel of the Feynman-Kac semigroup, 
is given in Section \ref{sec:HKE}.

In this paper, we adopt  the following notations. We use ``:=" as a way of definition. For two positive functions $f$ and $g$, notation $f\asymp g$ means that there is a constant $c\geq 1$ so that
$g/c \leq f \leq cg$, while notation $f\lesssim g$ (respectively, $f\gtrsim g$) means there is a constant $c>0$
so that $f\leq c g$ (respectively, $f\geq c g$). 

\section{Preliminaries and Main Result}
\label{sec:preliminary}

Throughout the remainder of this paper, we assume that $X$ is  a Hunt process  on $\R^d$ with $d\geq 2$ 
 having a jointly continuous transition density function $p(t, x, y)$ with respect to the Lebesgue 
measure on $\R^d$ and that
the two-sided estimates \eqref{e:1.3} holds for $p(t, x, y)$ on $(0, 1]\times \R^d \times \R^d$. 
Since we are concerned with heat kernel estimates of \eqref{e:1.3} on fixed time intervals, 
it is desirable to rewrite the estimates in the following equivalent but more compact form. 
This equivalent form \eqref{hkeofX2} is given in \cite{CHXZ}. For reader's convenience, we give a proof here.

\begin{lemma} 	\label{L:2.1}
	Two-sided estimates  \eqref{e:1.3} for $p(t, x, y)$ on $(0, 1]\times \R^d \times \R^d$ is equivalent 
	to the following. 
	There exist constants $C\geq1$ and $\lambda\in (0, 1]$ such that for $0< t\leq 1$ and $x,y\in\R^d$, 
	\begin{equation}
	\label{hkeofX2}
	C^{-1}\left(\Gamma_{c_2}(t;x-y)+\eta(t;x-y)\right)\leq p(t,x,y)\leq C\left(\Gamma_{c_4}(t;x-y)+\eta(t;x-y)\right),
	\end{equation}
	where
	\begin{equation}
		\label{eqn:gammaeta}
		\Gamma_\lambda(t;x):=t^{-d/2}e^{-\lambda|x|^2/t} 
		\quad \hbox{ and } \quad \eta(t;x):=\frac{t}{(t^{1/2}+|x|)^{d+\alpha}}.
	\end{equation}
\end{lemma}

\proof Note that  $t^{-d/2}\leq t^{-d/\alpha}$ for $t\in(0,1]$, and 
$$
\frac{1}{2}(a\wedge b+a\wedge c)\leq a\wedge(b+c)\leq a\wedge b+a\wedge c \quad \hbox{for }a,b,c>0.
$$
Thus for $\lambda >0$ and $t\in (0, 1]$, 
\begin{eqnarray*}
&& \left(t^{-d/2}\wedge t^{-d/\alpha}\right)\wedge\left(t^{-d/2}{\rm Exp}\left(-\frac{\lambda r^2}{t}\right)+t^{-d/\alpha}\wedge \frac{t}{r^{d+\alpha}}\right) \\
&\asymp & t^{-d/2}{\rm Exp}\left(-\frac{\lambda r^2}{t}\right)+t^{-d/2}\wedge t^{-d/\alpha}\wedge \frac{t}{r^{d+\alpha}} \\
&\asymp & \begin{cases}
t^{-d/2} {\rm Exp}\left(-\frac{\lambda r^2}{t}\right) + t^{-d/2}\wedge \frac{t}{r^{d+\alpha}} 
\asymp t^{-d/2} {\rm Exp}\left(-\frac{\lambda r^2}{t}\right) + \frac{t}{(t^{1/2}+r)^{d+\alpha}}
&\hbox{if } r\geq t^{1/2} \, (\, \geq t^{1/\alpha} \, ), \smallskip \\
 t^{-d/2}+t^{-d/2} \wedge \frac{t}{r^{d+\alpha}}\asymp t^{-d/2}  
&\hbox{if } t^{1/\alpha}<r\leq t^{1/2},  \smallskip \\
 t^{-d/2} + t^{-d/2} \wedge t^{-d/\alpha}\asymp t^{-d/2} 
&\hbox{if } 0<r < t^{1/\alpha} \, (\, \leq t^{1/2} \, ) 
\end{cases} \\
&\asymp &  t^{-d/2} {\rm Exp}\left(-\frac{\lambda r^2}{t}\right) + \frac{t}{(t^{1/2}+r)^{d+\alpha}},
\end{eqnarray*}
where the last line is due to the fact that for $0<r\leq t^{1/2}$, 
$\frac{t}{(t^{1/2}+r)^{d+\alpha}} \leq t^{1-(d+\alpha)/2} \leq t^{-d/2}$. 
This establishes the lemma.  
\qed

\par 
We now introduce some Kato classes for signed measures and for functions $F$ used in non-local Feynman-Kac perturbation.
For a $\sigma$-finite signed measure $\mu$, we use $\mu^+$ and $\mu^-$ to denote its positive and negative part
in its Jordan decomposition, and its total variation measure is given by $|\mu|:= \mu^++\mu^-$.
For a signed measure $\mu$ on $\R^d$, using the notations in (\ref{eqn:gammaeta}), we define
\begin{equation} \label{e:2.3}
N_\mu^{\alpha,\lambda}:=\sup_{x\in\R^d}\int_0^t\int_{\R^d}\left(\Gamma_\lambda(t;x-y)+\eta(t;x-y)\right)|\mu|(dy)ds.
\end{equation}
For a function $F(x, y)$ defined on $\R^d\times \R^d$ that vanishes along the diagonal, 
we define 
\begin{equation} \label{e:2.4}
N_F^{\alpha,\lambda}(t):=\sup_{y\in\R^d}\int_0^t\int_{\R^d\times\R^d}\left(\Gamma_\lambda(t;y-z)+\eta(t;y-z)\right)\frac{|F(z,w)|+|F(w,z)|}{|z-w|^{d+\alpha}}dwdzds.
\end{equation}

\begin{defn}\label{D:2.2}
\begin{description}
\item{\rm (i)}  A signed measure $\mu$ on $\R^d$ is said to be in the Kato class $\K_{\alpha}$ if 
\ $\lim_{t\downarrow0}N_\mu^{\alpha,\lambda}(t)=0$ for some and hence for all $\lambda >0$.
A measurable function $f$ on $\R^d$ is said to be in Kato class $\K_\alpha$ if $|f(x)| \mu (dx) \in \K_\alpha$.

\item{\rm (ii)} A bounded measurable function $F$ on $\R^d \times \R^d$ vanishing on the diagonal, is said to be 
Kato class $\J_{\alpha}$ if $\lim_{t\downarrow0}N_F^{\alpha,\lambda}(t)=0$ for some and hence for all $\lambda >0$. 
\end{description}
\end{defn}

\medskip

Clearly, if $F, G\in \J_{\alpha}$ and $a\in \R$, then so are $aF$, $e^F-1$, $ F+G$ and $ FG$. 
By H\"older inequality, it is easy to check that $L^\infty (\R^d) + L^p (\R^d)\subset \K_\alpha$
for every $p>d/2$. 
For $\mu \in \K_\alpha$ and $F\in \J_\alpha$, we can define an additive function of $X$ by  
$$
A_t^{\mu,F}=A_t^\mu+\sum_{0<s\leq t}F(X_{s-},X_s),
$$
where $A_t^\mu$ is a continuous additive functional of $X$ having $\mu$ as its Revuz measure.
It is easy to check that $A_t^{\mu,F}$ is well defined and is of  finite variations on compact time intervals. We can then 
define the following non-local Feynman-Kac semigroup of $X$ by 
\begin{equation}
	\label{eqn:FKsemigroup}
	T_t^{\mu,F}f(x)=\Ex_x\Big[{\rm Exp}(A_t^{\mu,F})f(X_t)\Big],\hspace{5mm}t\geq0, x\in\R^d.
\end{equation}

\medskip

The goal  of this paper is to study the stability of heat kernel estimates under the above non-local Feynman-Kac transform.

\begin{thm}
	\label{thm:main}
	Suppose $X$ is a Hunt process on $\R^d$ that has a jointly continuous transition density function $p(t, x, y)$
	with respect to the Lebesgue measure 
	and that the two-sided heat kernel estimates \eqref{e:1.3} holds for $p(t, x, y)$ on $(0, 1]\times \R^d \times \R^d$. 
	Let  $\mu\in \K_{\alpha}$ and $F(x, y)$ be a measurable function so that $F_1:=e^F-1 \in \J_\alpha$.
	Then the non-local Feynman-Kac semigroup $(T_t^{\mu,F};t\geq0)$ has a jointly continuous kernel $q(t,x,y)$ so that
	$T_t^{\mu,F} f(x)=\int_{\R^d} q(t, x, y) f(y) dy$ for every bounded Borel measurable function $f$ on $\R^d$.
	Moreover, there exist positive constants $c_1 $  and $K$ that depend on $(d,\alpha,C, N_{\alpha,F}^{\alpha, c_4},\norm{F_1 }_\infty)$
	so that for any $t>0$ and $x,y\in\R^d$,
	$$  q(t,x,y)\leq c_1 \, e^{Kt}\left( \Gamma_{2c_4/3}(t;x-y)+\eta(t;x-y) \right). 
	$$
	Here $c_4$ is the constant in \eqref{e:1.3}.
	If in addition, $F\in \J_\alpha$, then there exist positive constants $c_2 $, $\lambda_1$  and $K_1$ that depend on $(d,\alpha,C, N_{\alpha,F_1}^{\alpha, c_4},\norm{F }_\infty)$
	so that for any $t>0$ and $x,y\in\R^d$,
	$$  q(t,x,y)\geq c_2 \, e^{-K_1t}\left( \Gamma_{\lambda_1}(t;x-y)+\eta(t;x-y) \right). 
	$$
\end{thm}

\section{3P inequalities}
\label{sec:3P}

In this section we will establish various  3P type inequalities, which are key ingredients in the proof of Theorem \ref{thm:main}. Lemma \ref{lemma:3P}, Lemma \ref{lemma:gen3Pforeta},  Lemma \ref{lemma:gen3PforGamma},  and Lemma \ref{lemma:gen3PforGammaeta} are dealing with $\Gamma_c(t;x-y)$ and $\eta(t;x-y)$ as defined in (\ref{eqn:gammaeta}). Theorem \ref{thm:3Pforp} and Theorem \ref{thm:gen3Pforp} are the main results of this section.
\begin{lemma}
	\label{lemma:3P}
	For $0<s<t$, and $x,y,z\in\R^d$,
	\begin{enumerate}[\rm (i)]
		\item There exists a constant $C_1=C_1(d,\alpha)$ such that 
		\begin{equation}
			\label{ineq:eta}
			\eta(t-s;x-z)\eta(s;z-y)\leq C_1\eta(t;x-y)(\eta(t-s;x-z)+\eta(s;z-y)).
		\end{equation}
		\item For $0<a<b$, there exists a constant $C_2=C_2(d,a,b)$ such that for any measure $\mu$ on $\R^d$, 
		\begin{align}
			\label{ineq:gamma}
			&\int_0^t\int_{\R^d}\Gamma_a(t-s;x-z)\Gamma_b(s;z-y)|\mu|(dz)ds\notag\\
			&\hspace{1cm}\leq C_2\Gamma_a(t;x-y)\sup_{x\in\R^d}\int_0^t\int_{\R^d}\Gamma_c(s;x-y)|\mu|(dy)ds,
		\end{align}
		where $c=(b-a)\wedge\frac{a}{2}$.
		\item There exists a constant $C_3=C_3(d,\alpha,a)$ such that for any measure $\mu$ on $\R^d$,
		\begin{align}
			\label{ineq:gammaeta}
			&\int_0^t\int_{\R^d}\Gamma_a(t-s;x-z)\eta(s;z-y)|\mu|(dz)ds\notag\\
			&\leq C_3\bigg(\Gamma_a(t;x-y)\sup_{x\in\R^d}\int_0^t\int_{\R^d}\eta(s;x-y)|\mu|(dy)ds\notag\\
			&\hspace{1cm}+\eta(t;x-y)\sup_{x\in\R^d}\int_0^t\int_{\R^d}\Gamma_a(s;x-y)|\mu|(dy)ds\bigg).
		\end{align}
	\end{enumerate}
\end{lemma}

\proof  (i)	  For $0<s<t$ and $x,y,z\in\R^d$, we have
	\begin{align*}
	&\frac{\eta(t-s;x-z)\eta(s;z-y)}{\eta(t;x-y)}\\
	&\hspace{5mm}=\frac{(t-s)s}{t}\bigg(\frac{t^{1/2}+|x-y|}{((t-s)^{1/2}+|x-z|)(s^{1/2}+|z-y|)}\bigg)^{d+\alpha}\\
	&\hspace{5mm}\leq((t-s)\wedge s)\bigg(\frac{((t-s)+s)^{1/2}+|x-z|+|z-y|}{((t-s)^{1/2}+|x-z|)(s^{1/2}+|z-y|)}\bigg)^{d+\alpha}\\
	&\hspace{5mm}\leq2^{d+\alpha}((t-s)\wedge s)\bigg(\frac{1}{((t-s)^{1/2}+|x-z|)^{d+\alpha}}+\frac{1}{(s^{1/2}+|z-y|)^{d+\alpha}}\bigg)\\
	&\hspace{5mm}\leq 2^{d+\alpha}(\eta(t-s;x-z)+\eta(s;z-y)).
	\end{align*}

	(ii)  The proof for this part  is similar to that for  Lemma 3.1 in \cite{Z}. We first write
	\begin{align*}
	    J(t,x,y):=&\int_0^t\int_{\R^d}\Gamma_a(t-s;x-z)\Gamma_b(s;z-y)|\mu|(dz)ds\\
		=&\int_0^{\rho t}\int_{\R^d}\Gamma_a(t-s;x-z)\Gamma_b(s;z-y)|\mu|(dz)ds\\
		&+\int_{\rho t}^{t}\int_{\R^d}\Gamma_a(t-s;x-z)\Gamma_b(s;z-y)|\mu|(dz)ds
	\end{align*}
	Applying the elementary inequality
	\begin{equation}
		\label{ineq:elem}
		\frac{|x-z|^2}{t-s}+\frac{|z-y|^2}{s}\geq \frac{|x-y|^2}{t}, \text{ for }0<s<t,
	\end{equation}
	one can obtain
	\begin{align*}
		&\int_0^{\rho t}\int_{\R^d}\Gamma_a(t-s;x-z)\Gamma_b(s;z-y)|\mu|(dz)ds\\
		&=\int_0^{\rho t}\int_{\R^d}(t-s)^{-d/2}s^{-d/2}{\rm Exp}\Big(-a\frac{|x-z|^2}{t-s}\Big){\rm Exp}\Big(-b\frac{|z-y|^2}{s}\Big)|\mu|(dz)ds\\
		&= \int_0^{\rho t}\int_{\R^d}(t-s)^{-d/2}s^{-d/2}{\rm Exp}\Big(-a\big(\frac{|x-z|^2}{t-s}+\frac{|z-y|^2}{s}\big)\Big){\rm Exp}\Big(-(b-a)\frac{|z-y|^2}{s}\Big)|\mu|(dz)ds\\
		&\leq \int_0^{\rho t}\int_{\R^d}(t-s)^{-d/2}s^{-d/2}{\rm Exp}\Big(-a\frac{|x-y|^2}{t}\Big){\rm Exp}\Big(-(b-a)\frac{|z-y|^2}{s}\Big)|\mu|(dz)ds\\
		&\leq (1-\rho)^{-d/2}\Gamma_a(t;x-y)\int_0^{\rho t}\int_{\R^d}\Gamma_{b-a}(s;z-y)|\mu|(dz)ds.
	\end{align*}
	For the other term, by defining $U=\{|z-y|\geq |x-y|(a/b)^{1/2}\}$, we have
	\begin{align*}
		&\int_{\rho t}^{t}\int_{\R^d}\Gamma_a(t-s;x-z)\Gamma_b(s;z-y)|\mu|(dz)ds\\
		&=\int_{\rho t}^{t}\int_U\Gamma_a(t-s;x-z)\Gamma_b(s;z-y)|\mu|(dz)ds+\int_{\rho t}^{t}\int_{U^c}\Gamma_a(t-s;x-z)\Gamma_b(s;z-y)|\mu|(dz)ds\\
		&\leq (\rho t)^{-d/2}{\rm Exp}\Big(-\frac{a|x-y|^2}{t}\Big)\int_{\rho t}^{t}\int_U(t-s)^{-d/2}{\rm Exp}\Big(-\frac{a|x-z|^2}{t-s}\Big)|\mu|(dz)ds\\
		&\hspace{5mm}+(\rho t)^{-d/2}\int_{\rho t}^{t}\int_{U^c}(t-s)^{-d/2}{\rm Exp}\Big(-\frac{a|x-z|^2}{t-s}\Big)|\mu|(dz)ds.
	\end{align*}
	On $U^c$, we have the inequality,
	$$|x-z|\geq |x-y|-|y-z|\geq |x-y|\big(1-(a/b)^{1/2}\big),$$
	thus,
	\begin{align*}
		&\int_{\rho t}^{t}\int_{\R^d}\Gamma_a(t-s;x-z)\Gamma_b(s;z-y)|\mu|(dz)ds\\
		&\leq \rho^{-d/2}\Gamma_a(t;x-y)\int_{\rho t}^{t}\Gamma_a(t-s;x-z)|\mu|(dz)ds\\
		&\hspace{5mm}+(\rho t)^{-d/2}\int_{\rho t}^{t}\int_{U^c}(t-s)^{-d/2}{\rm Exp}\Big(-\frac{a|x-z|^2}{2(t-s)}\Big){\rm Exp}\Big(-\frac{a\big(1-(a/b)^{1/2}\big)^2|x-y|^2}{2(t-s)}\Big)|\mu|(dz)ds\\
		&\leq \rho^{-d/2}\Gamma_a(t;x-y)\int_{\rho t}^{t}\Gamma_a(t-s;x-z)|\mu|(dz)ds\\
		&\hspace{5mm}+(\rho t)^{-d/2}{\rm Exp}\Big(-\frac{a\big(1-(a/b)^{1/2}\big)^2|x-y|^2}{2(1-\rho)t}\Big)\int_{\rho t}^{t}\int_{U^c}(t-s)^{-d/2}{\rm Exp}\Big(-\frac{a|x-z|^2}{2(t-s)}\Big)|\mu|(dz)ds,
	\end{align*}
	by selecting $\rho$ such that $2(1-\rho)=\big(1-(a/b)^{1/2}\big)^2$, we would achieve the estimate in (\ref{ineq:gamma}), with $c=(b-a)\wedge\frac{a}{2}$, and $C_2$ depends on $d,a,b$.

	(iii) For $0<s<t$, if $|x-y|\leq t^{1/2}$, we have
	\begin{align*}
		&\Gamma_a(t-s;x-z)\leq 2^{d/2}t^{-d/2}\leq 2^{d/2}e^a\Gamma_a(t;x-y),&\text{ for }s\in(0,t/2];\\
		&\eta(s;z-y)\leq 4^{d+\alpha}\eta(t;x-y),&\text{ for }s\in(t/2,t).
	\end{align*}
	If $|x-y|>t^{1/2}$, consider on $V=\{|y-z|\geq|x-y|/2 \}$, we would have $\eta(s;z-y)\leq 2^{d+\alpha}\eta(t;x-y)$ for all $0<s<t$. 
	\par
	On $V^c$, $|x-z|\geq |x-y|-|y-z|\geq |x-y|/2$, we would have $\Gamma_a(t-s;x-z)\leq \gamma\Gamma_a(t;x-y)$, where $\gamma$ depends on $a, d$. 
	\par
	The estimate (\ref{ineq:gammaeta}) directly follows from the above discussion.
 \qed 

\par
Recall the definition of  $N_{\mu}^{\alpha,\lambda}$ from \eqref{e:2.3}.
 We next derive an integral 3P type inequality for $p(t,x,y)$ in small time, by using two-sided heat kernel estimates in Lemma \ref{L:2.1}.
 For notational convenience, let $\lambda=c_4$, where $c_4$ is the positive constant in \eqref{e:1.3}.

\begin{thm}
	\label{thm:3Pforp}
	For any $\mu\in\K_{\alpha}$, and any $(t,x,y)\in(0,1]\times\R^d\times\R^d$,
	\begin{equation}
	\label{eqn:3Pmu}
		\int_0^t\int_{\R^d}p(t-s,x,z)p_{2\lambda/3}(s,z,y)|\mu|(dz)ds\leq M_1p_{2\lambda/3}(t,x,y)N_\mu^{\alpha,\lambda/3}(t),
	\end{equation}
	where $M_1$ depends on $d,\alpha, C,\lambda$, and $p_{2\lambda/3}(t,x,y):=\Gamma_{2\lambda/3}(t;x-y)+\eta(t;x-y)$..
\end{thm}

\proof  By Lemma \ref{L:2.1},  we have
$$
p(t,x,y)\leq C\left(\Gamma_{\lambda}(t;x-y)+\eta(t;x-y)\right) \quad \hbox{for } t\in(0,1].
$$
Thus  for $(t,x,y)\in(0,1]\times\R^d\times\R^d$,
\begin{align*}
	&\int_0^t\int_{\R^d}p(t-s,x,z)p_{2\lambda/3}(s,z,y)|\mu|(dz)ds\\
	&\leq C\bigg(\int_0^t\int_{\R^d}\Gamma_{\lambda}(t-s;x-z)\Gamma_{2\lambda/3}(s;z-y)|\mu|(dy)ds+\int_0^t\int_{\R^d}\eta(t-s;x-z)\eta(s;z-y)|\mu|(dz)ds\\
	&\hspace{7mm}+\int_0^t\int_{\R^d}\Gamma_{\lambda}(t-s;x-z)\eta(s;z-y)|\mu|(dz)ds+\int_0^t\int_{\R^d}\eta(t-s;x-z)\Gamma_{2\lambda/3}(s;z-y)|\mu|(dz)ds\bigg)
\end{align*}
Applying Lemma \ref{lemma:3P}, we would have
\begin{align*}
&\int_0^t\int_{\R^d}p(t-s,x,z)p_{2\lambda/3}(s,z,y)|\mu|(dz)ds\\
&\leq C\bigg(C_2\Gamma_{2\lambda/3}(t,x-y)\sup_{x\in\R^d}\int_0^t\int_{\R^d}\Gamma_{\lambda/3}(s;x-y)|\mu|(dy)ds\\
&\hspace{1cm}+2C_1\eta(t;x-y)\sup_{x\in\R^d}\int_0^t\int_{\R^d}\eta(s;x-y)|\mu|(dy)ds\\
&\hspace{1cm}+2C_3\Gamma_{2\lambda/3}(t;x-y)\sup_{x\in\R^d}\int_0^t\int_{\R^d}\eta(s;x-y)|\mu|(dy)ds\\
&\hspace{1cm}+2C_3\eta(t;x-y)\sup_{x\in\R^d}\int_0^t\int_{\R^d}\Gamma_{2\lambda/3}(s;x-y)|\mu|(dy)ds\bigg),
\end{align*}
where $C_1, C_2, C_3$ depend on $d,\alpha, C,\lambda$. Altogether, let $M_1=C(2C_1\vee C_2\vee 2C_3)$,
we have
\begin{align*}
	&\int_0^t\int_{\R^d}p(t-s,x,z)p_{2\lambda/3}(s,z,y)|\mu|(dz)ds\\
	&\leq M_1\left(\Gamma_{2\lambda/3}(t,x-y)+\eta(t;x-y)\right)\sup_{x\in\R^d}\int_0^t\int_{\R^d}(\Gamma_{\lambda/3}(s;x-y)+\eta(s;x-y))|\mu|(dy)ds\\
	&=M_1p_{2\lambda/3}(t,x,y)N_\mu^{\alpha,\lambda/3}(t).
\end{align*}
\qed

\par 
We will use the following notations: for any $(x,y)\in \R^d\times\R^d$,
\begin{align*}
	&V_{x,y}:=\left\{(z,w)\in\R^d\times\R^d: |x-y|\geq 4(|y-w|\wedge|x-z|) \right\};\\
	&U_{x,y}:=V^c_{x,y}.
\end{align*}
First, similar as the discussion in \cite{CKS} (see Theorem 2.7), we could have the generalized integral 3P inequality for $\eta(t;x-y)$.

\begin{lemma}
	\label{lemma:gen3Pforeta}
	There exists a constant $C_4=C_4(\alpha,d)$ such that for any non-negative bounded function $F(x,y)$ on $\R^d\times\R^d$, the followings are true for $(t,x,y)\in(0,\infty)\times\R^d\times\R^d$.
	\begin{enumerate}[\rm (i)]
		\item If $|x-y|\leq t^{1/2}$, then 
		\begin{align*}
			&\int_0^t\int_{\R^d\times\R^d}\eta(t-s;x-z)\eta(s;w-y)\frac{F(z,w)}{|z-w|^{d+\alpha}}dzdwds\\
			&\leq C_4\eta(t;x-y)\int_0^t\int_{\R^d\times\R^d}(\eta(s;x-z)+\eta(s;w-y))\frac{F(z,w)}{|z-w|^{d+\alpha}}dzdwds.
		\end{align*}
		\item If $|x-y|>t^{1/2}$, then
		\begin{align*}
			&\int_0^t\int_{U_{x,y}}\eta(t-s;x-z)\eta(s;w-y)\frac{F(z,w)}{|z-w|^{d+\alpha}}dzdwds\\
			&\leq C_4\eta(t;x-y)\int_0^t\int_{U_{x,y}}(\eta(s;x-z)+\eta(s;w-y))\frac{F(z,w)}{|z-w|^{d+\alpha}}dzdwds.
		\end{align*}
		\item If $|x-y|>t^{1/2}$, then 
		$$\int_0^t\int_{V_{x,y}}\eta(t-s;x-z)\eta(s;w-y)\frac{F(z,w)}{|z-w|^{d+\alpha}}dzdwds\leq C_4\norm{F}_{\infty}\eta(t;x-y),$$
		where $\norm{F}_{\infty}$ denotes the $L^{\infty}$-norm of $F$ on $\R^d\times\R^d$.
	\end{enumerate}
\end{lemma}
Now we proceed to get the generalized integral 3P inequality for $\Gamma_c(t;x-y)$.

\begin{lemma}
	\label{lemma:gen3PforGamma}
	For $0<a<b$, there exists a constant $C_5=C_5(a,b,d)$ such that for any non-negative bounded function $F(x,y)$ on $\R^d\times\R^d$, the followings are true for $(t,x,y)\in(0,\infty)\times\R^d\times\R^d$.
	\begin{enumerate}[\rm (i)]
		\item If $|x-y|\leq t^{1/2}$, then 
		\begin{align*}
			&\int_0^t\int_{\R^d\times\R^d}\Gamma_a(t-s;x-z)\Gamma_b(s;w-y)\frac{F(z,w)}{|z-w|^{d+\alpha}}dzdwds\\
			&\leq C_5\Gamma_b(t;x-y)\int_0^t\int_{\R^d\times\R^d}(\Gamma_a(s;x-z)+\Gamma_b(s;w-y))\frac{F(z,w)}{|z-w|^{d+\alpha}}dzdwds.
		\end{align*}
		\item If $|x-y|> t^{1/2}$, then 
		\begin{align*}
		&\int_0^t\int_{U_{x,y}}\Gamma_a(t-s;x-z)\Gamma_b(s;w-y)\frac{F(z,w)}{|z-w|^{d+\alpha}}dzdwds\\
		&\leq C_5\Gamma_a(t;x-y)\int_0^t\int_{U_{x,y}}(\Gamma_a(s;x-z)+\Gamma_b(s;w-y))\frac{F(z,w)}{|z-w|^{d+\alpha}}dzdwds.
		\end{align*}
		\item If $|x-y|> t^{1/2}$, then 
		$$\int_0^t\int_{V_{x,y}}\Gamma_a(t-s;x-z)\Gamma_b(s;w-y)\frac{F(z,w)}{|z-w|^{d+\alpha}}dzdwds\leq C_5\norm{F}_{\infty}\eta(t;x-y).$$
	\end{enumerate}
\end{lemma}
 
\proof  (i) If $|x-y|\leq t^{1/2}$, we have $\Gamma_a(t-s;x-z)\leq 2^{d/2}e^{b}\Gamma_b(t;x-y)$ when $s\in(0,t/2])$; $\Gamma_b(s;w-y)\leq 2^{d/2}e^{b}\Gamma_b(t;x-y)$ when $s\in(t/2,t)$. Then (i) follows naturally.
	
	(ii)  If $|x-y|>t^{1/2}$, we let 
	\begin{align*}
		&U_1:=\{(z,w)\in\R^d\times\R^d:|y-w|>4^{-1}|x-y|,|y-w|\geq |x-z| \};\\
		&U_2:=\{(z,w)\in\R^d\times\R^d:|x-z|>4^{-1}|x-y| \}.
	\end{align*}
	Note that $\Gamma_b(s;w-y)\leq \gamma_1\Gamma_b(t;x-y)$ on $U_1$ for $s\in(0,t)$; $\Gamma_a(t-s;x-z)\leq \gamma_2\Gamma_a(t;x-y)$ on $U_2$ for $s\in(0,t)$, where $\gamma_1:=\gamma_1(b,d)$ and $\gamma_2:=\gamma_2(a,d)$. Since $U_{x,y}=U_1\cup U_2$, (ii) follows directly.
	
	(iii)  On $V_{x,y}$, $|z-w|\geq 2^{-1}|x-y|$. Hence 
	\begin{eqnarray*}
		&&  \int_0^t\int_{V_{x,y}}\Gamma_a(t-s;x-z)\Gamma_b(s;w-y)\frac{F(z,w)}{|z-w|^{d+\alpha}}dzdwds\\
		&\leq & 2^{d+\alpha}|x-y|^{-(d+\alpha)}\norm{F}_\infty\int_0^t\int_{V_{x,y}}\Gamma_a(t-s;x-z)\Gamma_b(s;w-y)dzdwds\\
		&\lesssim & \frac1{t} \eta(t;x-y)\norm{F}_\infty\int_0^t \Big( \int_{\R^d}  \Gamma_a(t-s; z)  dz \Big) \Big( \int_{\R^d}
		\Gamma_b(s;w) dw \Big) ds \\
		&\lesssim &   \eta(t;x-y)\norm{F}_\infty.
	\end{eqnarray*}
 This completes the proof of the lemma. \qed
 
\par 
We next establish a generalized integral 3P inequality  involving both $\Gamma_c(t;x-y)$ and $\eta(t;x-y)$.

\begin{lemma}
	\label{lemma:gen3PforGammaeta}
	There exists a constant $C_6=C_6(c,\alpha,d)$ such that for any non-negative bounded function $F(x,y)$ on $\R^d\times\R^d$, the followings are true for $(t,x,y)\in(0,1]\times\R^d\times\R^d$.
	\begin{enumerate}[\rm (i)]
		\item If $|x-y|\leq t^{1/2}$, then
		\begin{align*}
			&\int_0^t\int_{\R^d\times\R^d}\Gamma_c(t-s;x-z)\eta(s;w-y)\frac{F(z,w)}{|z-w|^{d+\alpha}}dzdwds\\
			&\leq C_6\bigg(\Gamma_c(t;x-y)\int_0^t\int_{\R^d\times\R^d}\eta(s;w-y)\frac{F(z,w)}{|z-w|^{d+\alpha}}dzdwds\\
			&\hspace{1cm}+\eta(t;x-y)\int_0^t\int_{\R^d\times\R^d}\Gamma_c(s;x-z)\frac{F(z,w)}{|z-w|^{d+\alpha}}dzdwds\bigg).
		\end{align*}
		\item If $|x-y|>t^{1/2}$, then
		\begin{align*}
			&\int_0^t\int_{U_{x,y}}\Gamma_c(t-s;x-z)\eta(s;w-y)\frac{F(z,w)}{|z-w|^{d+\alpha}}dzdwds\\
			&\leq C_6\bigg(\Gamma_c(t;x-y)\int_0^t\int_{U_{x,y}}\eta(s;w-y)\frac{F(z,w)}{|z-w|^{d+\alpha}}dzdwds\\
			&\hspace{1cm}+\eta(t;x-y)\int_0^t\int_{U_{x,y}}\Gamma_c(s;x-z)\frac{F(z,w)}{|z-w|^{d+\alpha}}dzdwds\bigg).
		\end{align*}
		\item If $|x-y|>t^{1/2}$, then
		$$\int_0^t\int_{V_{x,y}}\Gamma_c(t-s;x-z)\eta(s;w-y)\frac{F(z,w)}{|z-w|^{d+\alpha}}dzdwds\leq C_6\norm{F}_\infty\eta(t;x-y).$$
	\end{enumerate}
\end{lemma}

\proof (i) If $|x-y|\leq t^{1/2}$, we have $\Gamma_c(t-s;x-z)\leq 2^{d/2}e^c\Gamma_c(t;x-y)$ when $s\in(0,t/2]$; $\eta(s;w-y)\leq 4^{d+\alpha}\eta(t;x-y)$ when $s\in(t/2,t]$. Thus, we have (i) hold naturally.

	(ii)  If $|x-y|>t^{1/2}$, we continue to use the decomposition $U_{x,y}=U_1\cup U_2$ in the proof of Lemma \ref{lemma:gen3PforGamma}, and observe that $\Gamma_c(t-s;x-z)\leq \gamma_4(t;x-y)$ on $U_2$ for $s\in(0,t)$, where $\gamma_4$ depends on $d,c$. Thus, we first have
	\begin{align*}
		&\int_0^t\int_{U_{2}}\Gamma_c(t-s;x-z)\eta(s;w-y)\frac{F(z,w)}{|z-w|^{d+\alpha}}dzdwds\\
		&\leq\gamma_4\Gamma_c(t;x-y)\int_0^t\int_{U_{2}}\eta(s;w-y)\frac{F(z,w)}{|z-w|^{d+\alpha}}dzdwds.
	\end{align*}
	Also, observe that $\eta(s;w-y)\leq 4^{d+\alpha}\eta(t;x-y)$ when $s\in(0,t)$ and $(z,w)\in U_1$. Thus,
	\begin{align*}
	&\int_0^t\int_{U_{1}}\Gamma_c(t-s;x-z)\eta(s;w-y)\frac{F(z,w)}{|z-w|^{d+\alpha}}dzdwds\\
	&\leq 4^{d+\alpha}\eta(t;x-y)\int_0^t\int_{U_{1}}\Gamma_c(s;x-z)\frac{F(z,w)}{|z-w|^{d+\alpha}}dzdwds.
	\end{align*}
	Altogether, (ii) holds directly.

	(iii)  Note that on $V_{x,y}$, $|z-w|\geq 2^{-1}|x-y|$. Also, there exists $\gamma_5$ depending on $d,c,\alpha$ such that $\int_{\R^d}\Gamma_c(t;x-y)dy\leq \gamma_5$ and $\int_{\R^d}\eta(t;x-y)dy\leq \gamma_5t^{(2-\alpha)/2}$. Thus,
	\begin{align*}
		&\int_0^t\int_{V_{x,y}}\Gamma_c(t-s;x-z)\eta(s;w-y)\frac{F(z,w)}{|z-w|^{d+\alpha}}dzdwds\\
		&\leq 2^{d+\alpha}|x-y|^{-(d+\alpha)}\norm{F}_\infty\int_0^t\gamma_5^2s^{(2-\alpha)/2}ds\\
		&\leq 2^{d+\alpha}\gamma_5^2\norm{F}_\infty|x-y|^{-(d+\alpha)}t^{(4-\alpha)/2},
	\end{align*}
	for $t\leq 1$ and $|x-y|>t^{1/2}$, there exists $C_6$ depending on $d,\alpha,c$ such that
	$$\int_0^t\int_{V_{x,y}}\Gamma_c(t-s;x-z)\eta(s;w-y)\frac{F(z,w)}{|z-w|^{d+\alpha}}dzdwds\leq C_6\norm{F}_\infty\eta(t;x-y).$$
 \qed

\par
Recall that the definition of 	$N_F^{\alpha,\lambda}(t)$ from \eqref{e:2.4}.  
Note that a Hunt process $X_t$ admits a L\'evy system $(N(x, dy), H_t$), where $N(x, dy)$ is a kernel and $H_t$ is a positive continuous
additive functional of $X_t$;  that is, 
for any $x\in\R^d$, any stopping time $T$ and any non-negative measurable function $f$ on $[0,\infty)\times\R^d\times\R^d$, vanishing on the diagonal,
\begin{equation}
	\label{eqn:Lsys}
	\Ex_x\left[\sum_{s\leq T}f(s,X_{s-},X_s)\right]=\Ex_x\left[\int_0^T\int_{\R^d}f(s,X_s,y) N(X_s, dy) dH_s\right].
\end{equation}
Since $X_t$ has transition density function $p(t, x, y)$ with respect to the Lebesgue measure, it follows that
the Revuz measure $\mu_H$ of $H$ is absolutely continuous with respect to the Lebesgue measure.
So we can take $\mu_H (dx)=dx$, in other words, we can take $H_t\equiv t$. 
By two-sided heat kernel estimates \eqref{e:1.3} for the Hunt process $X$
and  the fact that $N(x, dy)$ is the weak limit of $p(t, x, y)dy/t$ as $t\to 0$,
we have
\begin{equation}\label{e:3.7} 
H_t=t \quad \hbox{ and } \quad N(x, dy) = \frac{c(x, y)}{|x-y|^{d+\alpha}} dy 
\end{equation}
for some measurable function $c(x, y)$ on $\R^d\times \R^d$ that is bounded between two positive constants.

\begin{thm} 
	\label{thm:gen3Pforp}
	Suppose $F (x, y)$ is a measurable function so that $F_1=e^F-1\in\J_{\alpha}$. There is a constant $M_2>0$ so that 
	 for any $(t,x,y)\in(0,1]\times\R^d\times\R^d$,
	 \begin{align}
	 \label{eqn:gen3Pforp}
		&\int_0^t\int_{\R^d\times\R^d}p(t-s,x,z)p_{2\lambda/3}(s,w,y)\frac{ |F_1|(z,w)}{|z-w|^{d+\alpha}}dzdwds\notag\\
		&\leq M_2\, p_{2\lambda/3}(t,x,y)\left(N_{F_1}^{\alpha,\lambda/3}(t)+\norm{F_1}_\infty\1_{\{|x-y|>t^{1/2}\}}\right),
		\end{align}
		 In particular, on $U_{x,y}=\{(z,w)\in\R^d\times\R^d:|x-y|\geq 4(|y-w|\wedge|x-z|)\}^c$,
		\begin{equation}
		\label{eqn:gen3PforponU}
			 \int_0^t\int_{U_{x,y}}p(t-s,x,z)p_{2\lambda/3}(s,w,y)\frac{ |F_1|(z,w)}{|z-w|^{d+\alpha}}dzdwds 
			 \leq M_2 \, p_{2\lambda/3}(t,x,y)N_{F_1}^{\alpha,\lambda/3}(t). 
		\end{equation} 
\end{thm}

\proof By Lemma \ref{L:2.1},
\begin{align*}
	&\int_0^t\int_{\R^d\times\R^d}p(t-s,x,z)p_{2\lambda/3}(s,w,y)\frac{ |F_1|(z,w)}{|z-w|^{d+\alpha}}dzdwds\\
	&\leq C\bigg(\int_0^t\int_{\R^d\times\R^d}\Gamma_{\lambda}(t-s;x-z)\Gamma_{2\lambda/3}(s;w-y)\frac{ |F_1|(z,w)}{|z-w|^{d+\alpha}}dzdwds\\
	&\hspace{1cm}+\int_0^t\int_{\R^d\times\R^d}\eta(t-s;x-z)\eta(s;w-y)\frac{ |F_1|(z,w)}{|z-w|^{d+\alpha}}dzdwds\\
	&\hspace{1cm}+\int_0^t\int_{\R^d\times\R^d}\Gamma_{\lambda}(t-s;x-z)\eta(s;w-y)\frac{ |F_1|(z,w)}{|z-w|^{d+\alpha}}dzdwds\\
	&\hspace{1cm}+\int_0^t\int_{\R^d\times\R^d}\eta(t-s;x-z)\Gamma_{2\lambda/3}(s;w-y)\frac{ |F_1|(z,w)}{|z-w|^{d+\alpha}}dzdwds\bigg).
\end{align*}
Applying (i) and (ii) in Lemma \ref{lemma:gen3Pforeta}, \ref{lemma:gen3PforGamma} and \ref{lemma:gen3PforGammaeta}, we first have for $|x-y|\leq t^{1/2}$, and for $\{|x-y|>t^{1/2}\}\cap U_{x,y}$,
$$
	 \int_0^t\int_{U_{x,y}}p(t-s,x,z)p_{2\lambda/3}(s,w,y)\frac{ |F_1|(z,w)}{|z-w|^{d+\alpha}}dzdwds 
	 \lesssim p_{2\lambda/3}(t,x,y)N_{F_1}^{\alpha,\lambda/3}(t).
$$
This establishes \eqref{eqn:gen3PforponU}.  For $|x-y|>t^{1/2}$ and $(z,w)\in V_{x,y}$, we apply (iii) in Lemma \ref{lemma:gen3Pforeta}, \ref{lemma:gen3PforGamma} and \ref{lemma:gen3PforGammaeta} to deduce 
$$ \int_0^t\int_{V_{x,y}}p(t-s,x,z)p_{2\lambda/3}(s,w,y)\frac{ |F_1|(z,w)}{|z-w|^{d+\alpha}}dzdwds
\lesssim \eta(t;x-y)\norm{F_1}_\infty.
$$
Hence inequality \eqref{eqn:gen3Pforp} holds.   \qed

\begin{lemma}\label{L:3.7} There is a constant $C>0$ so that for every $t\in (0, 2)$ and $x, y\in \R^d$, 
\begin{equation}
	\label{eqn:3Pforpc4}
	\int_{\R^d}p_{2\lambda/3}(t/2,x,z)
	p_{2\lambda/3}(t/2,z,y)dz\leq C \, p_{2\lambda/3}(t,x,y).
\end{equation}
\end{lemma}

 \proof It follows  from the 3P inequality for $\eta$ in Lemma \ref{lemma:3P}, we have
\begin{eqnarray*}
	&&\int_{\R^d}p_{2\lambda/3}(t/2,x,z)
	p_{2\lambda/3}(t/2,z,y)dz\\
	&\leq& \int_{\R^d}\Gamma_{2\lambda/3}(t/2;x-z)
	\Gamma_{2\lambda/3}(t/2;z-y)dz+\int_{\R^d}\eta(t/2;x-z)\eta(t/2;z-y)dz\\
	&&\quad +\int_{\R^d}\Gamma_{2\lambda/3}(t/2;x-z)
	\eta(t/2;z-y)dz+\int_{\R^d}\Gamma_{2\lambda/3}(t/2;z-y)
	\eta(t/2;x-z)dz\\
	&\lesssim  & \Gamma_{2\lambda/3}(t;x-y)+ \eta(t;x-y) +\int_{\R^d}\Gamma_{2\lambda/3}(t/2;x-z)
	\eta(t/2;z-y)dz \\
	&&\quad +\int_{\R^d}\Gamma_{2\lambda/3}(t/2;z-y)
	\eta(t/2;x-z)dz,
\end{eqnarray*}
for the second to the last term, when $|x-z|\geq \sqrt 2|x-y|/2$, $\Gamma_{2\lambda/3}(t/2;x-z)\leq 2^{d/2}\Gamma_{2\lambda/3}(t;x-y)$; when $|x-z|<\sqrt 2|x-y|/2$, then $|y-z|\geq |x-y|-|x-z|\geq \left(1-\frac{\sqrt 2}{2}\right)|x-y|$, we have $\eta(t/2;z-y)\leq \left(2/(2-\sqrt 2)\right)^{d+\alpha}\eta(t;x-y)$. Thus  
$$
\int_{\R^d}\Gamma_{2\lambda/3}(t/2;x-z) \eta(t/2;z-y)dz\lesssim p_{2\lambda/3}(t,x,y).
$$
With similar discussion for the last term, we   conclude that  
\eqref{eqn:3Pforpc4} holds. \qed

\section{Heat Kernel Estimates}
\label{sec:HKE}

In  the study of non-local Feynman-Kac perturbation, it is convenient to use Stieltjes exponential
rather than the standard exponential.
Recall that if $K_t$ is a right continuous function with left limits on $\R_+$ with
$K_0=1$ and $\Delta K_t:=K_t-K_{t-}>-1$ for every $t>0$, and if
$K_t$ is of finite variation on each compact time interval, then the
Stieltjes exponential ${\rm Exp} (K)_t$ of $K_t$ is the unique
solution $Z_t$ of
$$
 Z_t=1+ \int_{(0, t]} Z_{s-} dK_s, \quad t> 0.
$$
It is known that 
$$ {\rm Exp} (K)_t = e^{K^c_t} \prod_{0<s\leq t} (1+\Delta K_s),
$$
where $K_t^c$ denotes
the continuous part of $K_t$.
 The above formula gives a one-to-one correspondence between Stieltjes exponential and the natural exponential. 
The reason of ${\rm Exp} (K)_t$ being called the {\it Stieltjes}
exponential of $K_t$ is that,  by  \cite[p. 184]{DD},
${\rm Exp} (K)_t$ can be expressed as the following infinite sum of Lebesgue-Stieltjes integrals: 
\begin{equation}\label{e:4.1}
 {\rm Exp} (K)_t =1+ \sum_{n=1}^\infty
\int_{[0, t]} dK_{s_n} \int_{[0, s_n)} dK_{s_{n-1}} \cdots \int_{[0, s_2)} dK_{s_1}.
\end{equation}
The advantage of using the Stieltjes exponential ${\rm Exp} (K)_t$
over the usual exponential ${\rm Exp} (K_t)$  is the identity
\eqref{e:4.1}, which allows one to apply
the Markov property of $X$.

\subsection{Upper bound estimate}

Throughout this subsection, $\mu\in\K_{\alpha}$ and $F$ is a measurable function so that  $F_1:=e^F-1\in\J_{\alpha}$.
We will adopt the approach of \cite{CKS} to construct and derive its upper bound estimate  for the heat kernel 
of the non-local Feynman-Kac semigroup.
Define 
$$
N_{\mu,F_1}^{\alpha,\lambda}(t):=N_\mu^{\alpha,\lambda}(t)+N_{F_1}^{\alpha,\lambda}(t) 
$$
and let
\begin{equation}
	\label{eqn:K}
	K_t:=A_t^\mu+\sum_{s\leq t}F_1(X_{s-},X_s) . 
\end{equation}
  Then $\exp (A^\mu_t+\sum_{s\leq t} F_1(X_{s-},X_s)) = {\rm Exp} (K)_t$. So
  it follows from \eqref{e:4.1} that 
\begin{equation}
\label{eqn:semigroup}
T_t^{\mu,F}f(x)=P_tf(x)+\Ex_x\bigg[f(X_t)\sum_{n=1}^\infty\int_{[0,t]}dK_{s_n}\int_{[0,s_n)}dK_{s_{n-1}}\cdots\int_{[0,s_2)}dK_{s_1}\bigg] . 
\end{equation}
In view of  Theorem \ref{thm:3Pforp} and Theorem \ref{thm:gen3Pforp}, we can interchange the order of the expectation
and the unfinite sum (see the proof of Theorem \ref{thm:upperestimate} for details).  
Using the Markov property of $X$ and setting $h_1(s):=1$, $h_{n-1}(s) :=\int_{[0,s)}dK_{s_{n-1}}\cdots\int_{[0,s_2)}dK_{s_1}$,  we have
\begin{align}
\label{eqn:Tmu}
T_t^{\mu,F}f(x)
=&P_tf(x)+\sum_{n=1}^\infty\Ex_x\bigg[f(X_t)\int_{[0,t]}dK_{s_n}\int_{[0,s_n)}dK_{s_{n-1}}\cdots\int_{[0,s_2)}dK_{s_1}\bigg]\notag\\
=&P_tf(x)+\sum_{n=1}^\infty\Ex_x\bigg[\int_{[0,t]}P_{t-s_n}f(X_{s_n})dK_{s_n}\int_{[0,s_n)}dK_{s_{n-1}}\cdots\int_{[0,s_2)}dK_{s_1}\bigg]\notag\\
=&P_tf(x)+\sum_{n=1}^\infty\Ex_x\bigg[\int_{[0,t]}\bigg(\int_{[0,s_n)}P_{t-s_n}f(X_{s_n})h_{n-1}(s_{n-1})dK_{s_{n-1}}\bigg)dK_{s_n}\bigg]\notag\\
=&P_tf(x)+\Ex_x\bigg[\int_{[0,t)}\bigg(\Ex_{X_{s_{n-1}}}\Big[\int_{(0,t-s_{n-1}]}P_{t-s_{n-1}-r}f(X_r)dK_r\Big]\notag\\
&\hspace{35mm}\times\int_{[0,s_{n-1})}dK_{s_{n-2}}\cdots\int_{[0,s_2)}dK_{s_1}\bigg)dK_{s_{n-1}}\bigg] .
\end{align}
For any bounded measurable $g\geq 0$ on $[0,\infty)\times\R^d\times\R^d$,   by \Levy system of $X$ in  \eqref{eqn:Lsys}-\eqref{e:3.7}, 
\begin{align}
\label{eqn:gX}
\Ex_x\bigg[\int_{(0,s]}g(s-r,X_r)dK_r\bigg]=&\Ex_x\bigg[\int_{(0,s]}g(s-r,X_r)dA_r^\mu+\sum_{r\leq s}g(s-r,X_r)F_1(X_{r-},X_r)\bigg]\notag\\
=&\int_0^s\int_{\R^d}p(r,x,y)g(s-r,y)\mu(dy)dr\notag\\
&+\Ex_x\bigg[\int_0^s\bigg(\int_{\R^d}g(s-r,y)F_1(X_r,y)\frac{c(X_r,y)}{|X_r-y|^{d+\alpha}}dy\bigg)dr\bigg]\notag\\
=&\int_0^s\int_{\R^d}p(r,x,y)g(s-r,y)\mu(dy)dr\notag\\
&+\int_0^s\int_{\R^d}p(r,x,y)\bigg(\int_{\R^d}g(s-r,y)F_1(z,y)\frac{c(z,y)}{|z-y|^{d+\alpha}}dy\bigg)dzdr .
\end{align}
Define $p^{(0)}(t,x,y):=p(t,x,y)$, and for $k\geq1$,
\begin{align}
\label{eqn:pk}
p^{(k)}(t,x,y):=&\int_0^t\bigg(\int_{\R^d}p(t-s,x,z)p^{(k-1)}(s,z,y)\mu(dz)\bigg)ds\notag\\
&+\int_0^t\bigg(\int_{\R^d\times\R^d}p(t-s,x,z)p^{(k-1)}(s,w,y)\frac{c(z,w)F_1(z,w)}{|z-w|^{d+\alpha}}dzdw\bigg)ds.
\end{align}
Let 
\begin{equation}
q(t, x, y):= \sum_{n=0}^\infty p^{(k)} (t, x, y),
\end{equation}
which  will be shown in the proof of Theorem\ref{thm:upperestimate}  to be absolutely convergent under the assumption of $\mu \in \K_\alpha$
and $F_1\in \J_\alpha$. Then it follows from \eqref{eqn:gX}
and \eqref{eqn:gX} that
\begin{equation}
T_t^{\mu,F}f(x) = \int_{\R^d} q(t, x, y) f(y) dy.
\end{equation}
So  $q(t,x, y)$ is the heat kernel for the Feynman-Kac semigroup $\{T_t^{\mu,F}; t\geq 0\}$.
We will derive  upper bound estimate on $q(t, x, y)$ by estimating each $p^{(k)} (t, x, y)$.

\begin{lemma}
	\label{lemma:intofpk} There are  constants $C_0\geq 1$ and  $M\geq 1$ such that 
	for every  $k\geq 0$ and $(t,x)\in(0,1]\times\R^d$,  
	\begin{equation}
	\label{eqn:integralofpkonspace}
		\int_{\R^d}|p^{(k)}(t,x,y)|dy\leq C_0\left(MN_{\mu,F_1}^{\alpha,\lambda}(t)\right)^k.
	\end{equation}
\end{lemma}

\proof     We prove this lemma by induction. When $k=0$, by Lemma \ref{L:2.1},
we have the inequality hold naturally. Suppose (\ref{eqn:integralofpkonspace}) is true for $k-1$. Then by   (\ref{eqn:pk}),
\begin{eqnarray*}
	\int_{\R^k}|p^{(k)}(t,x,y)|dy
	&\leq&\int_{0}^{t}\bigg(\int_{\R^d}p(t-s,x,z)\Big(\int_{\R^d}p^{(k-1)}(s,z,y)dy\Big)|\mu|(dz)\bigg)ds\\
	&&+\int_{0}^{t}\bigg(\int_{\R^d\times\R^d}p(t-s,x,z)\frac{c(z,w)|F_1|(z,w)}{|z-w|^{d+\alpha}}
	  \Big(\int_{\R^d}p^{(k-1)}(s,w,y)dy\Big)dzdw\bigg)ds\\
	&\leq& C_0\left(MN_{\mu,F_1}^{\alpha,\lambda}(t)\right)^{k-1}\int_0^t\int_{\R^d}p(t-s,x,z)|\mu|(dz)ds\\
	&&+C_0\left(MN_{\mu,F_1}^{\alpha,\lambda}(t)\right)^{k-1}\int_0^t\int_{\R^d\times\R^d}p(t-s,x,z)
	 \frac{c(z,w)|F_1|(z,w)}{|z-w|^{d+\alpha}}dzdwds\\
	&\leq&C_0  C (1+  \| c\|_\infty) M^{k-1}\left(N_{\mu,F_1}^{\alpha,\lambda}(t)\right)^k\leq C_0\left(MN_{\mu,F_1}^{\alpha,\lambda}(t)\right)^k, 
\end{eqnarray*}
if we increase the value of $M$ if necessary so that $M\geq C (1+ \|c\|_\infty)$. Here $C\geq 1$ is the constant in Lemma \ref{L:2.1}. The lemma is proved. 
\qed

\begin{lemma}
	\label{lemma:pk}
	For any $k\geq 0$ and $(t,x,y)\in(0,1]\times\R^d\times\R^d$,  
\begin{equation}\label{e:4.10}
		|p^{(k)}(t,x,y)|\leq Cp_{2\lambda/3}(t,x,y)\left((MN_{\mu,F_1}^{\alpha,\lambda/3}(t))^k+k\norm{F_1}_\infty M(MN_{\mu,F_1}^{\alpha,\lambda/3}(t))^{k-1}\right) ,
\end{equation}
where $C\geq 1 $ and $M\geq 1$ are the constants in Lemma \ref{L:2.1} and Lemma \ref{lemma:intofpk}, respectively.
\end{lemma}

\proof Inequality holds trivially for $k=0$. Suppose it is true for $k-1\geq 0$, then if $|x-y|\leq t^{1/2}$, using the induction hypothesis and applying Theorem \ref{thm:3Pforp} and Theorem \ref{thm:gen3Pforp},
\begin{eqnarray*}
	|p^{(k)}(t,x,y)|
	&\leq&\int_0^t\bigg(\int_{\R^d}p(t-s,x,z)|p^{(k-1)}(s,z,y)||\mu|(dz)\bigg)ds\\
	&& +\int_0^t\bigg(\int_{\R^d\times\R^d}p(t-s,x,z)\frac{c(z,w)|F_1|(z,w)}{|z-w|^{d+\alpha}}|p^{(k-1)}(s,w,y)|dzdw\bigg)ds\\
	&\leq&C\left((MN_{\mu,F_1}^{\alpha,\lambda/3})^{k-1}+(k-1)\norm{F_1}_\infty M(MN_{\mu,F_1}^{\alpha,\lambda/3})^{k-2}\right)\\
	&& \times\bigg(\int_0^t\bigg(\int_{\R^d}p(t-s,x,z)p_{2\lambda/3}(s,z,y)|\mu|(dz)\bigg)ds\\
	&& \hspace{1cm}+\int_0^t\bigg(\int_{\R^d\times\R^d}p(t-s,x,z)p_{2\lambda/3}(s,w,y)
	     \frac{c(z,w)|F_1|(z,w)}{|z-w|^{d+\alpha}}dzdw\bigg)ds\bigg)\\
	&\leq&Cp_{2\lambda/3}(t,x,y)\left((MN_{\mu,F_1}^{\alpha,\lambda/3})^{k-1}+(k-1)\norm{F_1}_\infty M(MN_{\mu,F_1}^{\alpha,\lambda/3})^{k-2}\right)MN_{\mu,F_1}^{\alpha,\lambda/3}(t).
\end{eqnarray*}
If $|x-y|>t^{1/2}$, we have
\begin{eqnarray*}
	|p^{(k)}(t,x,y)| &\leq&\int_0^t\bigg(\int_{\R^d}p(t-s,x,z)|p^{(k-1)}(s,z,y)||\mu|(dz)\bigg)ds\\
	&&+\int_0^t\bigg(\int_{U_{x,y}}p(t-s,x,z)\frac{c(z,w)|F_1|(z,w)}{|z-w|^{d+\alpha}}|p^{(k-1)}(s,w,y)|dzdw\bigg)ds\\
	&&+\int_0^t\bigg(\int_{V_{x,y}}p(t-s,x,z)\frac{c(z,w)|F_1|(z,w)}{|z-w|^{d+\alpha}}|p^{(k-1)}(s,w,y)|dzdw\bigg)ds\\
	&=&J_1+J_2+J_3 .
\end{eqnarray*}
Applying Theorem \ref{thm:3Pforp} to $J_1$ and Theorem \ref{thm:gen3Pforp} to $J_2$,
\begin{eqnarray*}
	J_1+J_2
	&\leq&C\left((MN_{\mu,F_1}^{\alpha,\lambda/3})^{k-1}+(k-1)\norm{F_1}_\infty M(MN_{\mu,F_1}^{\alpha,\lambda/3})^{k-2}\right)\\
	&& \quad \times\bigg( \int_0^t\int_{\R^d}p(t-s,x,z)p_{2\lambda/3}(t,x,y)|\mu|(dz)ds\\
	&&\qquad  \quad +\int_0^t\bigg(\int_{U_{x,y}}p(t-s,x,z)p_{2\lambda/3}(s,w,y)\frac{c(z,w)|F_1|(z,w)}{|z-w|^{d+\alpha}}dzdw\bigg)ds \bigg)\\
	&\leq& Cp_{2\lambda/3}(t,x,y)\left((MN_{\mu,F_1}^{\alpha,\lambda/3})^{k-1}+(k-1)\norm{F_1}_\infty M(MN_{\mu,F_1}^{\alpha,\lambda/3})^{k-2}\right)MN_{\mu,F_1}^{\alpha,\lambda/3} . 
\end{eqnarray*}
 For $J_3$, use the fact that $|z-w|\geq 2^{-1}|x-y|$ and Lemma \ref{lemma:intofpk},
\begin{eqnarray*}
	J_3 &\leq&\frac{2^{d+\alpha}\norm{F_1}_\infty}{|x-y|^{d+\alpha}}\int_0^t\bigg(\int_{\R^d\times\R^d}p(t-s,x,z)|p^{(k-1)}(s,w,y)|dzdw\bigg)ds\\
	&\leq&2^{d+\alpha}\norm{F_1}_\infty\frac{t}{|x-y|^{d+\alpha}}C_0^2(MN_{\mu,F_1}^{\alpha,\lambda})^{k-1}\\
	&\leq& M\norm{F_1}_{\infty}p_{2\lambda/3}(t,x,y)(MN_{\mu,F_1}^{\alpha,\lambda/3})^{k-1} . 
\end{eqnarray*}
This  completes the proof.
\qed

\par 
The following result gives the existence and the desired upper bound estimates of the heat kernel for the non-local Feynman-Kac semigroup
$\{ T_t^{\mu,F};t\geq 0\}$, as stated in  Theorem \ref{thm:main}.

\begin{thm}
	\label{thm:upperestimate}
	The series $\sum_{k=0}^{\infty}p^{(k)}(t,x,y)$ converges absolutely to a jointly continuous function $q(t,x,y)$ on $(0,\infty)\times\R^d\times\R^d$. The function
	$q(t,x,y)$ is the integral kernel (or, heat kernel) for the Feynman-Kac semigroup $\{T_t^{\mu,F};  t\geq 0\}$, 
	and there exists a constant $K$ depending on $d,\alpha, \norm{F_1}_\infty$ and the constants $ C$ and $\lambda:=c_4$ in Lemma \ref{L:2.1} such that
	\begin{equation}
		q(t,x,y)\leq e^{Kt}p_{2\lambda/3}(t,x,y)  \quad \hbox{for every } (t,x,y)\in(0,\infty)\times\R^d\times\R^d.
	\end{equation}
 \end{thm}

\proof  Let  $\wh p^{(k)}(t,x,y)$ be defined as in  \eqref{eqn:pk}  
 but with $|\mu|$ and  $|F_1|$ in place of $\mu$ and $F_1$; that is, $\wh p^{(0)}(t, x, y)= p(t, x, y)$, 
 and for $k\geq 1$, 
\begin{align}
	\label{eqn:hatpk}
	\wh p^{(k)}(t,x,y):=&\int_0^t\bigg(\int_{\R^d}  p(t-s,x,z)\wh p^{(k-1)}(s,z,y)|\mu|(dz)\bigg)ds\notag\\
	&+\int_0^t\bigg(\int_{\R^d\times\R^d}  p(t-s,x,z)\wh p^{(k-1)}(s,w,y)\frac{c(z,w)|F_1|(z,w)}{|z-w|^{d+\alpha}}dzdw\bigg)ds.
\end{align}
Clearly, $|p^{(k)}(t,x,y)|\leq\wh p^{(k)}(t,x,y)$ and by the proof of   Lemma \ref{lemma:pk}, there is a constant $0<t_1\leq 1$ so that 
 such that  $N_{\mu,F_1}^{\alpha,\lambda/3}(t_1)\leq (2M)^{-1}$ and that 
\begin{eqnarray}
\label{eqn:hatq}
	\wh q(t,x,y)&:=& \sum_{k=0}^{\infty}\wh p^{(k)}(t,x,y)  \nonumber \\
	 & \leq& Cp_{2\lambda/3}(t,x,y)+Cp_{2\lambda/3}(t,x,y)\sum_{k=1}^{\infty}\left((MN_{\mu,F_1}^{\alpha,\lambda/3}(t))^k+k\norm{F_1}_\infty M(MN_{\mu,F_1}^{\alpha,\lambda/3}(t))^{k-1}\right)\notag\\
	&\leq&Cp_{2\lambda/3}(t,x,y)+Cp_{2\lambda/3}(t,x,y)\left(1+4\norm{F_1}_\infty M\right)\notag\\
	&\leq&C(2+4\norm{F_1}_{\infty}M)p_{2\lambda/3}(t,x,y)=:\gamma_1 p_{2\lambda/3}(t,x,y).
\end{eqnarray}
This in particular implies that $\wh q(t,x,y)$ is  jointly continuous  on $(0, t_1]\times \R^d \times \R^d$.  Repeating the procedure (\ref{eqn:semigroup}), (\ref{eqn:Tmu}) and (\ref{eqn:gX}) with $|\mu|$, $|F_1|$ in place of $\mu$, $F_1$.  and by Fubini's theorem, we have for any bounded function $f\geq 0$ on $\R^d$
and $t\in (0, t_1]$,
\begin{equation}
\label{semigroupofhatp}
	T_t f(x):= \Ex_x\Big [f(X_t)\text{Exp}\Big(A^{|\mu|}+\sum_{s\leq\cdot}|F_1|(X_{s-},X_s)\Big)_t\Big]=\int_{\R^d}\wh q(t,x,y)f(y)dy.
\end{equation}
Note that $T_{t} \circ T_s = T_{t+s}$ for any $t, s \geq 0$. 
 Extend the definition of $\wh q(t, x, y)$ to $(0, 2t_1]\times \R^d\times \R^d$ by 
 $$
 \wh q (t+s, x, y) =\int_{\R^d} \wh q(t, x, z) \wh q (s, z, y) dz
 $$
 for $s, t \in (0, t_1]$. The above is well defined and, in view of  
 \eqref{eqn:hatq} and Lemma \ref{L:3.7}, $\wh q(t,x,y)$ is jointly continuous on $[0,2t_1]\times\R^d\times\R^d$ and 
 there is  constant $\gamma_2$ so that $\wh p(t, x, y) \leq \gamma_2 p_{2\lambda/3}(t,x,y) $ on
 $(0, 2t_1]\times \R^d \times \R^d$.  Clearly,
 $$ T_t f (x)= \int_{\R^d} \wh q (t,x, y) f(y) dy
 $$
 for every $f\geq 0$ on $\R^d$ and $(t, x)\in (0, 2t_1]\times \R^d$. 
 Repeat the above procedure, we can extend $\wh q(t, x, y)$ to be a jointly continuous function on 
  $[0,\infty)\times\R^d\times\R^d$ so that \eqref{semigroupofhatp} holds for  every $f\geq 0$ on $\R^d$ and $(t, x)\in (0, \infty) \times \R^d$,
  and that    there exists a constant $K>0 $ depending on $d,\alpha,C,\lambda,C_0,\norm{F_1}_\infty,M$ so that for any $t>0$ and $x,y\in\R^d$
$$\wh q(t,x,y)\leq e^{Kt}p_{2\lambda/3}(t,x,y).$$ 
This proves the theorem as $q(t, x, y)\leq \wh q(t, x, y)$. 
  \qed
 
\subsection{Lower bound estimate}

In this subsection, we assume $\mu \in \K_\alpha$ and    $F\in\J_{\alpha}$. Clearly, $F_1 :=e^F -1\in \J_\alpha$. 
Due to the presence of the Gaussian component in \eqref{e:1.3}, the approach in \cite{CKS} of obtaining lower bound estimates
for $q(t, x, y)$ is not applicable here. We will employ a probabilistic approach from \cite{CW, CZh} to get the desired lower bound estimates. 

Let  $\wt p^{(1)}(t,x,y)$ be defined as in \eqref{eqn:hatpk} but with $|F|$ in place of $|F_1|$. 
Thus by \eqref{eqn:hatq}, there is a constant $\gamma >0$ so that  $$\wt{p}^{(1)} (t,x,y)\leq \gamma p_{2\lambda/3}(t,x,y) 
\quad \hbox{for } (t,x, y) \in (0, 1] \times \R^d \times \R^d.
$$
In particular,  there is a constant $K_1>0$ so that $ \wt{p}^1(t,x,y)\leq K_1t^{-d/2}$  for $t\in(0,1]$ and $|x-y|\leq \sqrt t$.
 On the other hand, it follows from \eqref{hkeofX2} that  there exists a constant 
 $ \wt C\geq 1$ so that  
 $$
 \wt C^{-1}t^{-d/2}\leq p(t,x,y)\leq\wt Ct^{-d/2} \quad \hbox{for  every } t\in(0,1 \hbox{  and } |x-y|\leq \sqrt {t}.
 $$ 
 Let $k\geq 2$ be an integer so that $k\geq 2K_1\wt C$. Then for every  $t\in(0,1]$ and $x, y\in \R^d$ with $|x-y|\leq \sqrt t$, 
\begin{equation}
\label{eqn:phatp}
p(t,x,y)-\frac{1}{k}\wt {p}^1(t,x,y)\geq \frac{1}{2\wt C}t^{-d/2}\geq \frac{1}{2\wt C^2}p(t,x,y).
\end{equation}  
Note that 
$$\Ex_x\left[A_t^{|\mu|,|F|}f(X_t)\right] =\int_{\R^d}\wt p^{(1)}(t,x,y)f(y)dy. 
$$
Using the elementary inequality that  
$$
1-A_t^{|\mu|,|F|}/k\leq   \exp (- A_t^{|\mu|, |F|}/k) \leq \exp (A_t^{\mu,F}/k), 
$$
we have for  any  ball   $B(x,r)$ centered at $x$ with radius $r$ and any $(t, y) \in(0,1] \times \R^d$,  
$$
\frac{1}{|B(x,r)|}\Ex_y\left[(1-A_t^{|\mu|,|F|}/k)\1_{B(x,r)}(X_t)\right]\leq
\frac{1}{ | B(x,r)|}\Ex_y\left[ \exp(A_t^{\mu,F}/k)\1_{B(x,r)}(X_t)\right].
$$
Hence by  \eqref{eqn:phatp}  and   H\"older's inequality, we have for $0<t\leq 1$ and $|x-y|\leq \sqrt t$, 
\begin{align*}
	&\frac{1}{2\wt C^2}\frac{1}{B(x,r)}\Ex_y[\1_{B(x,r)}(X_t)]\leq \frac{1}{B(x,r)}\Ex_y[ \exp (A_t^{\mu,F}/k)\1_{B(x,r)}(X_t)]\\
	&\leq  \left(\frac{1}{B(x,r)}\Ex_y[ \exp (A_t^{\mu,F})\1_{B(x,r)}(X_t)]\right)^{1/k}\left(\frac{1}{B(x,r)}\Ex_y[\1_{B(x,r)}(X_t)]\right)^{1-1/k}.
\end{align*}
Thus  
$$
\frac{1}{2^k\wt C^{2k}}\frac{1}{B(x,r)}\Ex_y[\1_{B(x,r)}(X_t)]\leq\frac{1}{B(x,r)}\Ex_y[{\rm Exp}(A_t^{\mu,F})\1_{B(x,r)}(X_t)].
$$
By taking  $r\downarrow0$,  we conclude from above as well as  Lemma \ref{L:2.1} 
that  
\begin{equation}\label{e:4.16a}
q(t,x,y)\geq 2^{-k}\wt C^{-2k} p(t, x, y) \gtrsim t^{-d/2} \qquad \hbox{for every } t\in (0, 1] \hbox{ and } |x-y|\leq\sqrt {t}.
\end{equation}
By a standard chaining argument (see, e.g., \cite{FS}),  it follows that  there exist constants $K_2,\lambda_1>0$ so that 
\begin{equation}
\label{e:4.14}
q(t,x,y)\geq K_2\Gamma_{\lambda_1}(t;x-y) \quad \hbox{for  } (t, x, y)\in (0, 1]\times \R^d \times \R^d.
\end{equation}  

To get the jumping component in the lower bound estimate for $q(t, x, y)$, we consider a sub-Markovian semigroup 
$\{Q_t; t\geq 0\}$ defined by  
$$
Q_tf(x):=\Ex_x\left[{\rm Exp}\left(-A_t^{|\mu|}-\sum_{s\leq t}|F|(X_{s-},X_s)\right)f(X_t)\right].
$$
Since $|\mu \in \K_\alpha$ and $|F|\in \J_\alpha$, we know that  
  $\{Q_t;t\geq 0\}$ has a jointly continuous transition kernel $\bar p(t,x, y)$.
  Clearly, $q(t, x, y) \geq \bar p(t,x, y) $ for every $t>0$ and $x, y\in \R^d$.
  Since $\{Q_t;t\geq 0\}$   forms  a Feller semigroup,    there exists a  Feller process $Y=\{Y_t,\P_x,x\in\R^d,\zeta^Y\}$ such that $Q_tf(x)=\Ex_x[f(Y_t)]$.
 We will derive a lower bound estimate on $q(t,x, y)$ through the Feller process $Y$. 

It follows from the definition of $\mu\in\K_{\alpha}$ and $F\in\J_{\alpha}$ that
$\sup_{x\in \R^d} \Ex_x \left[ A^{ |\mu|, |F|}_{ 1} \right]<\infty$. Thus by Jensen's inequality,  
\begin{equation}\label{e:4.16} 
   \inf_{x\in \R^d} \Ex_x \left[{\rm Exp}(A^{-|\mu|,-|F|})_{ 1} \right] 
\geq  {\rm Exp} \left(-\sup_{x\in \R^d} \Ex_x \left[A^{ |\mu|, |F|}_{1} \right] \right) =:\gamma_0>0.
\end{equation}
Let $\eta$ be the random time whose distribution is determined by 
$\P_x (\zeta >t) = \Ex_x \left[{\rm Exp}(A^{-|\mu|,-|F|})_{ t} \right] $.
We can couple the processes $X$ and $Y$  in such a way that on $\{\eta >t\}$, 
$Y_s=X_s$ for every $s\leq t$. 

We define   first hitting time and exit time of a  Borel set $D\subset\R^d$ by   $X_t$ and $ Y_t$ as follows; 
\begin{align*}
	&\sigma^X_D:=\inf\{s\geq 0,X_s\in D\},\hspace{1cm}\tau^X_D:=\inf\{s\geq0,X_s\notin D \};\\
	&\sigma^Y_D:=\inf\{s\geq 0,Y_s\in D\},\hspace{1cm}\tau^Y_D:=\inf\{s\geq0,Y_s\notin D \}.
\end{align*}

\begin{lemma}
	\label{L:4.4} 
	Let $\gamma_0\in (0, 1)$ be the constant in \eqref{e:4.16}. 
	There exists a constant $\kappa_0\in(0,1)$ depending on $d,C,\lambda,\alpha,R,\gamma_0$ such that for any $0<r\leq 1$,  
	\begin{equation}
	\label{eqn:exit}
		\sup_{x\in\R^d}\P_x \left(\tau^X_{B(x,r)}\leq \kappa_0 r^2 \right)\leq  {\gamma_0}/{2}.
	\end{equation}
	Consequently, for every $x\in \R^d$ and $r\in (0, 1]$,  
	\begin{equation} \label{e:4.18}
		 \P_x \left(   \tau^Y_{B(x,r)}> \kappa_0 r^2 \right)\geq  {\gamma_0}/{2}.
	\end{equation}
\end{lemma}

\proof  First note that by \eqref{e:1.3}, for every $s\leq t\leq 1$, $x\in \R^d$ and $r>0$,
$$ 
\int_{B(x, r/)^c} p(s, x, y) dy \leq c  \left(\int_{B(0, r/2\sqrt{s})^c} e^{-c_4|z|^2}dz + s r^{-\alpha} \right)
\leq c \left( e^{-c_4 r^2/8t} + tr^{-\alpha} \right).
$$
Let $\kappa_0>0$ sufficiently small so that $c \left( e^{-c_4  /8\kappa_0} + \kappa_0 \right) <\gamma_0/4$.
Then by taking $t=\kappa_0r^2 $, we have from the above that for every $x\in \R^d$ and $r\in (0, 1]$, 
\begin{equation}\label{e:4.19} 
 \sup_{s\in (0, \kappa_0 r^2] } \int_{B(x, r/2)^c} p(s, x, y) dy \leq \gamma_0/4.
\end{equation}

For simplicity,  denote $\tau^X_{B(x,r)}$ by $\tau$. We have by the strong Markov property of $X_t$ and 
\eqref{e:4.19} that
\begin{eqnarray*}
	\P_x(\tau\leq \kappa_0 r^2 ) &\leq&\P_x\left(\tau\leq \kappa_0 r^2;  \, X_{\kappa_0 r^2}\in B(x,r/2)\right)
	+\P_x\left(X_{\kappa_0 r^2}\notin B(x,r/2)\right)\\
	&\leq&\P_x\left(\P_{X_\tau}(|X_{\kappa_0 r^2-\tau} -X_0|\geq r/2); \, \tau\leq \kappa_0 r^2)\right)+\frac{\gamma_0}4\\
	&\leq& \frac{\gamma_0}4+ \frac{\gamma_0}4 = \frac{\gamma_0}2. 
\end{eqnarray*}
Hence
\begin{eqnarray*}
  \P_x \left(   \tau^Y_{B(x,r)}> \kappa_0 r^2 \right) 
  &\geq& \P_x \left(  \eta >  \kappa_0 r^2  \hbox{ and } \tau^X_{B(x,r)} > \kappa_0 r^2 \right) \\
 & \geq & \P_x \left(  \eta >  \kappa_0 r^2 \right) -\P_x \left( \tau^X_{B(x,r)} \leq \kappa_0 r^2 \right)
 \geq \gamma_0/2, 
\end{eqnarray*}
where in the last inequality, we used \eqref{e:4.16}. 
 \qed

 \begin{lemma}
	\label{L:4.5}
Let $0< \kappa_0 <1 $ be the constant in Lemma \ref{L:4.4}.
There exists a constant $\gamma_1 >0$ so that
for any $r>0$  and $x_0, y_0\in \R^d$ with
$|y_0-x_0|\geq 3r$,  
	\begin{equation}
		\label{eqn:hittingtime}
		\P_{x_0} \left(\sigma^Y_{B(y_0,r)}\leq \kappa_0 r^2 \right)\geq \gamma_1\frac{r^{d+2}}{|y_0-x_0|^{d+\alpha}}.
	\end{equation}
\end{lemma}

\proof    Define $f(x, y)= 1_{B(x_0, r)} (x) 1_{B(y_0, r)} (y)$. 
Then
$$ 
M_t:=\sum_{s\leq t\wedge (\kappa r^2)} f(X_{s-}, X_s) -\int_0^{ t\wedge (\kappa_0 r^2)}
 f(X_s, y) N(X_s, dy) ds, \quad t\geq 0,
$$
is a martingale additive functional of $X$ that is  uniformly integrable under $\P_x$ for every $x\in \R^d$.  Let
$$
A_t:= A^{-|\mu|, |F|}_t = -A^{|\mu|}_t -\sum_{s\leq t} F( X_{s-}, X_s), \quad t\geq 0,
$$
which is a non-increasing additive functional of $X$. By stochastic integration by parts formula, 
$$
e^{A_t} M_t  = \int_0^t e^{A_{s-}} dM_s + \int_0^t M_{s-} de^{A_s} + \sum_{s\leq t} (e^{A_s} - e^{A_{s-}}) (M_s-M_{s-}).
$$
For $\tau:=\tau^X_{B(x_0, r)} \wedge (\kappa_0 r^2)$, $M_t= -\int_0^{ t\wedge (\kappa_0 r^2)}
 f(X_s, y) N(X_s, dy) ds\leq 0$ for $t\in [0, \tau)$. It follows that
 $$ e^{A_\tau} M_\tau \geq \int_0^{\tau} e^{A_{s-}} dM_s  + (e^{A_\tau} -e^{A_{\tau-}}) M_{\tau}.
   $$
 Thus 
 $$ 
 \Ex_{x_0} \left[ e^{A_{\tau-}} M_\tau \right] \geq \Ex_{x_0} \int_0^{\tau} e^{A_{s-}} dM_s =0.
 $$
 This together with \eqref{e:3.7} implies that 
 \begin{eqnarray*}
 \Ex_{x_0} \left[ e^{A_{\tau-}} 1_{B(y_0, r)} (X_\tau ) \right] 
 &\geq & \Ex_{x_0} \left[ e^{A_{\tau-}} \int_0^{ \tau\wedge   (\kappa_0 r^2)}
 \int_{B(y_0, r)} N(X_s, dy) ds \right] \\
 &\geq &  \frac{c \,\kappa_0 r^2\, r^d}{|x_0-y_0|^{d+\alpha}} \Ex_{x_0} \left[ e^{A_{\kappa_0 r^2 }};  \tau^X_{B(x_0, r)} \geq \kappa_0 r^2     \right] \\
 &=& \frac{c \,\kappa_0 r^{d+2} }{|x_0-y_0|^{d+\alpha}} \P_{x_0} \left( \tau_{B(x_0, r)}^Y\geq \kappa_0 r^2  \right) \\
 &\geq & \frac{c \, \gamma_0 \kappa_0 r^{d+2} }{2|x_0-y_0|^{d+\alpha}},
 \end{eqnarray*}
where the last inequality is due to \eqref{e:4.18}. Consequently, 
\begin{eqnarray*}
\P_{x_0} \left( \sigma^Y_{B(y_0,r)}\leq \kappa_0 r^2 \right)
&\geq & \P_{x_0} \left( \tau^Y_{B(x_0, r)} \leq \kappa_0 r^2  \hbox{ and } Y_{\tau^Y_{B(x_0, r)} } \in  B(y_0,r)  \right) \\
&=& \Ex_{x_0} \left[ e^{A_{\tau }} 1_{B(y_0, r)} (X_\tau ) \right]  \\
&\geq &  e^{-\| F\|_\infty} \, \Ex_{x_0} \left[ e^{A_{\tau-}} 1_{B(y_0, r)} (X_\tau ) \right]  \\
&\geq &  \frac{c \, e^{-\| F\|_\infty} \, \gamma_0 \kappa_0 r^{d+2} }{2|x_0-y_0|^{d+\alpha}}.  
\end{eqnarray*}
The lemma is proved.  \qed

We now derive  lower bound heat kernel estimate for the heat kernel $q(t, x, y)$ of the  Feynman-Kac semigroup $\{ T^{\mu, F}_t; t\geq 0\}$. 

\begin{thm}\label{T:4.6} 
	Suppose $\mu\in\K_{\alpha}$ and $F$ is a measurable function in $\J_{\alpha}$. 
	Then  there exist positive constants $\wt K\geq 1$ and $\lambda_1>0$ 
	depending on $(d,\alpha, \lambda,N_{\alpha,F}^{\alpha,\lambda/3},\norm{F}_\infty )$ 
	and the constants in \eqref{e:1.3} such that
	\begin{equation}\label{e:4.23}
		\wt K^{-1} \, p_{\lambda_1}(t,x,y)\leq q(t,x,y)\leq \wt K \, p_{2\lambda/3}(t,x,y) 
	\end{equation}
	for $(t,x,y)\in(0,1]\times\R^d\times\R^d$.
\end{thm}

\proof  The upper bound estimates follows from Theorem \ref{thm:upperestimate} so it remains to establish the lower bound estimate
for $q(t, x, y)$.
If $|x-y|\leq\sqrt t$, the desired lower bound heat kernel estimate follows from \eqref{e:4.14}. 
So it suffices to consider the case that $|x-y|>\sqrt t$.  Set $r=\sqrt{ t}/3$.
It follows from  Lemma \ref{L:4.4} and Lemma \ref{L:4.5} that  
\begin{eqnarray*}
	\P_x \left( Y_{2\kappa_0r^2}\in B(y,2r) \right)
	&\geq&\P_x\Big(\sigma^Y_{B(y,r)}<\kappa_0r^2;\sup_{s\in[\sigma,\sigma+\kappa_0r^2]}|Y_s-Y_\sigma|<r \Big)\\
	&=&\Ex_x \Big[\P_{Y_\sigma}\Big(\sup_{s\in[\sigma,\sigma+\kappa_0r^2]}|Y_s-Y_\sigma|<r\Big);\sigma<\kappa_0r^2\Big]\\
	&\geq&\P_x \left( \tau^Y_{B(x,r)}>\kappa_0r^2 \right)\P_x \left(\sigma^Y_{B(y,r)}<\kappa_0r^2 \right)\\
	&\geq&\frac{\gamma_0\gamma_1}{2}\frac{r^{d+2}}{|x-y|^{d+\alpha}}.
\end{eqnarray*}
 Thus
  $$
  \int_{B(y,2r)}q(2\kappa_0r^2,x,z)dz \geq \P_x \left(Y_{2\kappa_0r^2}\in B(y,2r) \right)
  \geq\frac{\gamma_0\gamma_1r^{d+2}}{2|x-y|^{d+\alpha}}.
 $$
Since $|y-z|<2r<\sqrt{t-2\kappa_0r^2}$, one has by (\ref{e:4.16a}) that    
\begin{eqnarray*}
	q(t,x,y) &\geq&\int_{B(y,2r)}q(2\kappa_0r^2,x,z)q(t-2\kappa_0r^2,z,y)dz\\
	&\geq&\inf_{z\in B(y,2r)}q(t-2\kappa_0r^2,y,z)\frac{\gamma_0\gamma_1r^{d+2}}{2|x-y|^{d+\alpha}}\\
	&\geq&K_2e^{-\lambda_1}t^{-d/2}\frac{\gamma_0\gamma_1r^{d+2}}{2|x-y|^{d+\alpha}} \\
	&\geq & K_2e^{-\lambda_1}\frac{\gamma_0\gamma_1}{2\cdot3^{d+2}}\eta(t;x-y).
\end{eqnarray*}
This together with \eqref{e:4.14} establishes the lower bound estimate for $q(t, x, y)$ in \eqref{e:4.23}.  \qed


\vskip 0.3truein 

Department of Mathematics, University of Washington, Seattle,
WA 98195, USA

Email: \texttt{zqchen@uw.edu}

Email: \texttt{lidanw@uw.edu}

\end{document}